\documentstyle[12pt]{article}
\newtheorem{propose}{Proposition}
\newtheorem{lemma}[propose]{Lemma}
\renewcommand{\P}{\mbox{$\bf P$}}   
\newcommand{\Q}{\mbox{$\bf Q$}}   
\newcommand{\C}{\mbox{$\bf C$}}     
\newcommand{\Z}{\mbox{$\bf Z$}}     
\newcommand{\Pic}{{\rm Pic}\,}      
\newcommand{\Ext}{{\rm Ext}\,}      
\newcommand{\supp}{{\rm supp}\,}
\newcommand{\sext}{\mbox{${\cal E}xt\,$}}  
\newcommand{\into}{\hookrightarrow}
\newenvironment{proof}{
                       \trivlist \item[\hskip \labelsep{\bf Proof}:]
                      }{
                        \hfill$\Box$\endtrivlist
                      }
\newenvironment{rmk}{
                        \trivlist \item[\hskip \labelsep{\bf Remark}:]
                      }{
                        \endtrivlist}

\newenvironment{claim}{
                        \trivlist \item[\hskip \labelsep{\bf Claim}:]
                      }{
                        \endtrivlist}
\newcounter{pts}

\def\thepts{\roman{pts}}
\def\labelpts{(\thepts)}
\newcommand{\eg}{{\it e.g.\/},\  }
\newcommand{\ie}{{\it i.e.\/},\  }
\newcommand{\rank}{{\rm rank}\,}    
\def\cF{{\cal F}}
\def\cI{{\cal I}}
\def\cL{{\cal L}}
\def\cM{{\cal M}}
\def\cO{{\cal O}}
\def\rarr{\rightarrow}
\newcommand{\pee}{\P^{n+1}}
\newcounter{try}
\author{N. Mohan Kumar\\
        Department of Mathematics\\
   	Washington University\\
        St. Louis, MO, 63130\\
          kumar@math.wustl.edu}
\title{Construction of rank two vector bundles on $\P^4$ in positive
characteristic}
\begin{document}
\maketitle

\section{Introduction}

Vector bundles on Projective spaces have been the subject of many
papers and many  results in this direction are known. For a
somewhat dated account, the reader may see \cite{OSS}. One of the most
interesting problems in this area is the study of small rank bundles
on Projective spaces. For example, a conjecture by R.~Hartshorne
\cite{rHartshorne} states that there are no small rank vector bundles
on Projective spaces, other than direct sum of line bundles.  The
solution to this tantalising problem still seems remote, though very
many results are known. Let me restrict my attention to rank 2 bundles
for the moment.  Many interesting bundles of rank 2 are known over
Projective spaces of dimension 2 and 3. But over $\P^4$, essentially
the only interesting one known is the well known Horrocks-Mumford
bundle \cite{HorMum}. There are also some intersting ones in
characteristic 2, discovered by Tango \cite{Tango} and G.~Horrocks
\cite{gHorrocks}.

In this paper, we shall deal with this problem and prove a criterion
relating bundles on $\P^{n+1}$ to bundles on $\P^n$. This condition
on certain bundles over $\P^n$
is necessary and sufficient for the existence of bundles on
$\P^{n+1}$. Though this criterion is not very pleasant, it allows you
to restrict your attention to bundles just on $\P^n$ to construct
bundles on $\P^{n+1}$. Using this criterion (which has nothing to do
with the characteristic of the field), we construct many rank two
bundles on $\P^4$ over a field of positive characteristic. (For the
Chern classes of these bundles, see section \ref{chern}).
Unfortunately we have not been able to extend our construction to
complex numbers, though I feel it should be possible. The
construction follows closely what we did in \cite{MK}.

\section{Genral remarks}\label{remarks}

Let me start with a word about notation. We will have to deal with
maps, $\phi: M\to M\otimes L$ often in this article, where  $M$ is a
sheaf and $L$, a line bundle. Then it makes sense to talk about
$\phi\otimes {\rm Id\/}: M\otimes L'\to M\otimes L\otimes L'$ for any
line bundle $L'$. We will denote this map also for brevity by $\phi$.
Also it makes sense to talk about $\phi^i: M\to M\otimes L^i$, by
composing $\phi$. So we will talk about $\phi$ as an {\em
endomorphism}, though strictly speaking, it is not. It also make sense
to say when such a map is {\em nilpotent}, by saying that $\phi^n=0$
for some $n$.

We make the following transparent observation:
\begin{rmk}
If $\pi: X\to Y$ is a finite map then the category of sheaves on $Y$
which are $\pi_*\cO_X$ modules is the same as the category of shevaes
of the form $\pi_*\cF$ where $\cF$ is a sheaf on $X$ (with the
appropriate homomorphisms).
\end{rmk}

A typical case where we plan to apply this is when $X$ is the
$m^{th}$ order thickening of $\P^n\subset\P^{n+1}$ and $Y=\P^n$ with
$\pi$ the projection from a point in $\P^{n+1}$ away from this
hyperplane. We had applied this in \cite{MK} in a similar but
slightly different context.

So let $X$ be the $m^{th}$ order thickening of $\P^n\subset\P^{n+1}$
and $Y=\P^n$ and $\pi: X\to Y$ the projection from a point away from
this hyperplane.

\begin{enumerate}\label{obv1}
\item Then $\pi_*\cO_X=\oplus_{i=0}^{m-1} \cO_Y(-i)$. Thus we see that
giving a sheaf on $X$ is equivalent to giving a sheaf $\cF$ on $Y$ and
an endomorphism $\phi: \cF\to\cF(1)$ with $\phi^m=0$.
\item Let $(E_1,\phi_1)$ and $(E_2,\phi_2)$ be two such sheaves on
$Y$. Giving a map $\psi:E_1\to E_2$ which commutes with the
$\phi_i$'s  is
equivalent to giving a map between the corresponding sheaves on $X$.
\item We want to apply this to the special case when $E_1$ arises as
the restriction of a direct sum of line bundles say $F$ on
$\P^{n+1}$. We will write for clarity $G$, for $F$ restricted to
$\P^n$. Then we see that
$$E_1=\pi_*(F\mid_X)=\oplus_{i=0}^{m-1} G(-i)$$ as before. $\phi_1$ can
be identified with the map which just shifts the blocks. That is to
say the $\phi_1$ takes $G(-i)$ to the corresponding $G(-i)$ in
$E_1(1)$ as identitiy. (Of course $G(-m+1)$ goes to zero).

\item Thus giving a map from $F$ as above to the sheaf corresponding
to $(E,\phi)$ is just giving a map $\theta: G\to E$. Because then we
get for any $i$ a map $G(-i)\to E$ by taking the induced map
$G(-i)\to E(-i)$ and then composing it with $\phi^i$.
\end{enumerate}

Now let $E$ be any vector bundle on $\pee$ and let $Y=\P^n$ be any
hyperplane. Then by Quillen--Suslin Theorem [\eg see\cite{Q}], E
restricted to the complement of this hyperplane is free. Thus, if we
denote by $F=\cO_{\pee}^r$ where $r=\rank E$, then we have an exact
sequence,
$$0\to E(-i)\to F\to \cF\to 0$$
for some integer $i$ and $\cF$ is a coherent sheaf on some $X$ as
above ($m^{th}$ order thickening of the hyperplane for some $m$).
Since we will be primarily interested in deciding when such an $E$ is
a direct sum of line bundles, we may as well rename $E(-i)$ by $E$.
Using $\pi$ as above we get a coherent sheaf $\pi_*\cF=M$ on $Y$. But
$\cF$ has homological dimension one and therefore by
Auslander-Buchsbaum Theorem [\eg see \cite{Mat}], $M$ is a vector
bundle. As above, we also have a nilpotent map $\phi:M\to M(1)$.
Letting $G=\cO_Y^r$ we also have a map $\psi: G\to M$ since the
distinguished $r$ sections of $\cF$ give $r$ sections of $M$. Further
by the surjectivity of the above map from $F\to\cF$, we see that
$$\phi(M(-1))+\psi(G)=M.$$ What we want to state is the converse of
this remark.

\begin{lemma}\label{extension}
Let $Y=\P^n$, $G$ be a rank $r$ bundle on $Y$ which is a direct sum
of line bundles. Assume $M$ is a vector bundle on $Y$ with
a nilpotent map $\phi: M\to M(1)$ and a map $\psi: G\to M$ such that
$\phi(M(-1))+\psi(G)=M$. Then there exists a vector bundle $E$ of
rank $r$ on $\pee$ and an exact sequence,
\begin{equation}\label{rem1}
0\to E\to F\to \cF\to 0
\end{equation}
where $F$ is direct sum of $r$ line bundles on $\pee$ with
$F\mid_{\P^n}=G$ and $\pi_*\cF=M$.
\end{lemma}
\begin{proof}
Proof is obvious using the remarks above.
\end{proof}

\begin{rmk}

The above lemma can be also thought of as a criterion for extending
vector bundles from $\P^n$ to $\pee$. In other words, given a vector
bundle $E$ of rank $r$ on $\P^n$, it can be extended to $\pee$ as a
vector bundle if and only if there exists a vector bundle $M$ over
$\P^n$ with a nilpotent endomorphism $\phi: M\to M(1)$, a map $\psi:
G\to M$ where $G$ is a direct sum of $r$ line bundles and an exact
sequence,

$$0\to E\to M(-1)\oplus G\stackrel{(\phi,\psi)}{\longrightarrow} M\to
0.$$

\end{rmk}

\begin{lemma}\label{nonsplit}
Let the notation be as in the above lemma and assume $n\geq 2$. If
$M$ is not a direct sum of line bundles then $E$ is not a direct sum
of line bundles. Conversely, if $r\leq n$, and $M$ is a direct sum of
line bundles, then so is $E$.
\end{lemma}

\begin{proof}

Assume $M$ is not a direct sum of line bundles. Then by Horrock's
criterion. [\eg see\cite{OSS}], $H^i(M(l))\neq 0$ for some $i$, with
$0<i<n$ and some $l\in\Z$.  $H^j(F(l))=0\quad \forall j, 0<j\leq n$,
since $F$ is a direct sum of line bundles. Therefore from our exact
sequence \ref{rem1}, $H^i(\cF(l))=H^{i+1}(E(l))$. Also
$H^i(\cF(l))=H^i(M(l))$ since $\pi$ is a finite map from $\supp \cF\to
Y$. Thus $H^{i+1}(E(l))\neq 0$ and since $0<i+1<n+1$, we see that $E$
is not a direct sum of line bundles.

Conversely, assume that $M$ is a direct sum of line bundles. Exactly
as before, we get $H^i(E(l))=0 \quad\forall l, 1<i\leq n$. By
duality, $H^i(E^*(l))=0\quad\forall l, 0<i<n$. Thus by the Syzygy
theorem \cite{GE} rank of $E^*=r>n$ or $E^*$ is a direct sum of
line bundles.
\end{proof}

\section{The Construction}

In this section, we will outline the construction of bundles $M$ on
$\P^3$ as described above. More generally, let $d\in\Z$ be any
integer. We will construct a bundle $M$ on $\P^3$ with a nilpotent
endomorphism $\phi: M\to M(d)$ and a map $\psi: G\to M$ where $G$ is
the direct sum of two line bundles such that $\phi(M(-d))+\psi(G)=M$
and $M$ is not a direct sum of line bundles, over a field of positive
characteristic.

\begin{rmk}
The only intersting cases are $d=-1,0,1$. If we have $M$'s for these
values, then by taking the pull back of these by finite maps
$\P^3\to\P^3$, we can construct bundles for all $d$. The case $d=0$
was treated in \cite{MK}.
\end{rmk}

Let $p>0$ be the characteristic of our algebraically closed field.
Choose positive numbers $N,k,l$ so that $N-k,N-l$ both positive,
$4pkl> d^2$ and
$p(k+l)=(p-1)N+d$. Let $x,y,z,t$ be the homogeneous co-ordinates of
$\P^3$. Let $A=x^kz^{N-k}+y^lt^{N-l}$. Let $C_i$ be the curve defined
by the vanishing of $x^{pk}, y^{pl}$ and $A^i$, for $1\leq i\leq p$.
Let $C$ be the curve (line) defined by $x=y=0$.

\begin{claim}

\begin{enumerate}
\item $C_i$'s are local complete intersection curves for $1\leq i\leq
p$ and $C_p$ is a complete intersection of $x^{pk}, y^{pl}$.

\item $\omega_{C_i}\cong \cO_{C_i}((i-1)N+d-4)$ for $1\leq i\leq p$
where $\omega$ as usual denote the dualising sheaf.

\item Thus by Serre Construction [\eg see\cite{OSS}], if we denote by
$L_i=\cO_Y(-(i-1)N-d)$ for all $i$, then we have exact sequences,
$$0\to L_i\stackrel{\alpha_i}{\longrightarrow}
M_i\stackrel{\beta_i}{\longrightarrow} \cI_{C_i}\to 0,$$
where $M_i$ are rank two vector bundles on $Y$ for $1\leq i\leq p$.
In fact we can arrange these extensions to fit into the following
commutative diagrams,
\[
\begin{array}{ccccccccc}

0&\to& L_i&{\stackrel{\alpha_i}{\longrightarrow}}&
M_i&{\stackrel{\beta_i}{\longrightarrow}}& \cI_{C_i}&\to& 0\\
&&\downarrow\cdot A&&\downarrow\eta_i&&\downarrow&&\\
0&\to& L_{i-1}&{\stackrel{\alpha_{i-1}}{\longrightarrow}}&
M_{i-1}&{\stackrel{\beta_{i-1}}{\longrightarrow}}& \cI_{C_{i-1}}&\to& 0
\end{array}
\]
where $L_i\to L_{i-1}$ is multiplication by $A$ and $\cI_{C_i}\to
\cI_{C_{i-1}}$ is the natural inclusion of ideals.

\item There exists a nilpotent endomorphism $\phi:M_1\to M_1(d)$
given as follows: Notice that $L_1=\cO_Y(-d)$. So we can identify
$\cI_{C_1}\subset L_1(d)$ and then define $\phi=\alpha_1\beta_1$.
\item $M_i/\eta_{i+1}(M_{i+1})$ is annihilated by $A$,
$1\leq i<p$. Thus the natural map $M_i\otimes\cO_Y(-N)\to M_i$ got by
multiplication by $A$ factors through $\eta_{i+1}$.

\item We have maps $g_i: L_{i+1}\to M_i(-d)$ for $1\leq i\leq p$ by
lifting $A^i$. \ie the composite $\beta_i\circ g_i$ is just given by
the element $A^i\in\cI_{C_i}$. We can arrange these maps so that
$\eta_i\circ g_i=g_{i-1}\circ A$ and $\phi\circ g_1=\alpha_1\circ A$.
\setcounter{try}{\theenumi}
\end{enumerate}
\end{claim}

The claim is proved exactly as in \cite{MK}. In fact the whole
construction is similar to that in \cite{MK}. So let me defer these
proofs and continue the construction assuming the above claim. Let
$$\cL=L_p\oplus L_{p-1}(-d)\oplus \cdots \oplus L_2(-d(p-2))$$
and
$$\cM=M_p\oplus M_{p-1}(-d)\oplus\cdots\oplus M_1(-d(p-1))$$
We have a map $f:\cL\to\cM$ given by sending
$(x_p,x_{p-1},\ldots,x_2)\in\cL$ to
$$(-\alpha_p(x_p),-\alpha_{p-1}(x_{p-1})+g_{p-1}(x_p),\ldots,
-\alpha_2(x_2)+g_2(x_3), g_1(x_2))\in\cM$$
Let the cokernel be called $M$.

\begin{claim}
\begin{enumerate}
\addtocounter{enumi}{\thetry}
\item $M$ is a rank $p+1$ vector bundle on $Y$.
\setcounter{try}{\theenumi}
\end{enumerate}
\end{claim}

We have an endomorphism $\theta: \cM\to \cM(d)$ given by,
$$(x_p,x_{p-1},\ldots, x_1)\mapsto (0, \eta_p(x_p),
\ldots,\eta_3(x_3),\eta_2(x_2)+ \phi(x_1))$$
Since $\phi^2=0$, one can easily see that $\theta^{p+1}=0$.

\begin{claim}
\begin{enumerate}
\addtocounter{enumi}{\thetry}
\item $\theta$ descends to a nilpotent endomorphism, $\varphi:M\to M(d)$
\setcounter{try}{\theenumi}
\end{enumerate}
\end{claim}

We have the natural map $\psi:M_p\to M$. Notice that since $C_p$ is a
complete intersection, $M_p$ is the direct sum of two line bundles.
Finally we have,

\begin{claim}
\begin{enumerate}
\addtocounter{enumi}{\thetry}
\item $$\psi(M_p)+\varphi(M(-d))=M$$
\item $M_1$ is not a direct sum of line bundles and hence neither
is $M$.
\setcounter{try}{\theenumi}
\end{enumerate}
\end{claim}

By taking $d=1$ in the above construction, we get a rank 2 bundle $E$
on $\P^4$ by lemma \ref{extension}. Since $M$ is not a direct sum of
line bundles, by lemma \ref{nonsplit}, $E$ is not a direct sum of line
bundles.

\subsection{Computation of Chern classes}\label{chern}
Our vector bundle $E$ on $\P^4$ is given by the exact sequence,
$$0\to E\to \cO(-pk)\oplus\cO(-pl)\to \cF\to 0$$
where $M=\pi_*\cF$ is the vector bundle on $\P^3$ as we have
constructed.
Notice that we are looking at the case $d=1$. So
to compute the Chern classes of $E$, we may as well restrict to a
general linear space of dimension 2, since rank of $E$ is 2. On this
$\P^2$ we will compute the class of $E$ in $K_0$. We have our
distinguished $\P^3\subset\P^4$ and the curve $C\subset\P^3$. So by
choosing our linear space generally, we may assume that it does not
meet this curve. Then we have $\P^1\subset\P^2$, after intersecting
with this linear space. We will denote by the same letters
restrictions of all our vector bundles. Since $\cI_{C_i}=\cO_{\P^1}$
now, we see that $[\cM_i]=[\cO]+[L_i]$ in $K_0(\P^1)$. Thus,
$$[M]=[\cO]+[\cO(-1)]+\ldots+[\cO(-(p-1))]+[L_1(-(p-1))]$$
$$=p[\cO]+[\cO(-{p(p+1)\over 2})].$$
Thus on $\P^2$, we see that,
$$[\cF]=p[\cO]-p[\cO(-1)]+[\cO(-{p(p+1)\over 2})]-[\cO(-1-{p(p+1)\over 2})].$$
So we get $[E]$ to be,
$$[\cO(-pk)]+[\cO(-pl)]-p[\cO]+p[\cO(-1)]-[\cO(-{p(p+1)\over 2})]+
[\cO(-1-{p(p+1) \over 2})].$$
Now an easy computation will show that,
$$c_1(E)=-1-p(k+l+1)$$
$$c_2(E)=p(p+1)(k+l)+p^2kl$$

For instance, taking $p=2$ and $k=l=1$, we get,
$$c_1(E)=-7, \quad c_2(E)=16.$$

By choosing appropriate $k,l$, one can construct vectorbundles with
$c_1^2>4c_2$ for example, in any characteristic, $p>0$. For instance,
let $s\geq 1$ be any integer and $k=1, l=ps-s$. Then $N=ps+1$ and the
corresponding rank two vector bundle has
$$c_1^2-4c_2=\alpha s^2+\beta s+\gamma$$
where $\alpha,\beta,\gamma$ depend only on $p$ and
$\alpha=p^2(p-1)^2>0$. Thus by choosing $s>>0$, we can make the
vector bundle to be of the required type.

\section{Proofs of the claims}
\begin{enumerate}

\item $C_p$ is a complete itersection is clear, since $A^p$ is in the
ideal generated by $x^{pk}$ and $y^{pl}$. To check the rest, we need
only look at points where either $z\neq 0$ or $t\neq 0$. If $z\neq 0$
one sees immediately that $x^{pk}\in (y^{pl}, A^i)$.

\item This is done by descending induction on $i$. For $i=p$, since
$C_p$ is a complete intersection of $x^{pk},y^{pl}$, this is obvious.
So assume result proved for all  $p\geq i>1$. Then we have an exact
seqence,
$$0\to \cO_{C_{i-1}}(-N)\to \cO_{C_i}\to \cO_{C_1}\to 0$$
which we dualise to get,
$$0\to \omega_{C_1}\to \omega_{C_i}\to \omega_{C_{i-1}}(N)\to 0$$
Since we already know from 1) that the last term is a line bundle on
$C_{i-1}$ and then the proof is clear.

\item This is just Serre construction.
Assume we have constructed the exact sequences upto $i-1$ with the
commutative diagrams, the first one is just by the usual Serre
construction. By taking the natural map $L_i\to L_{i-1}$ given by
multiplication by $A$, we get a map,
$$H^0(\cO_{C_i})=\Ext^1(\cI_{C_i}, L_i)\to
\Ext^1(\cI_{C_i}, L_{i-1})=H^0(\cO_{C_i}(N))$$
which is just multiplication by $A$. We also have a natural map,
induced from the inclusion, $\cI_{C_i}\subset \cI_{C_{i-1}}$,
$$H^0(\cO_{C_{i-1}})=\Ext^1(\cI_{C_{i-1}}, L_{i-1})\to
\Ext^1(\cI_{C_i}, L_{i-1})=H^0(\cO_{C_i}(N))$$
In this case also, it is clear that the element `1'$\in
H^0(\cO_{C_{i-1}})$ goes to `A'$\in H^0(\cO_{C_i}(N))$, which is also
the image of `1'$\in H^0(\cO_{C_i})$ by multiplication by $A$. But
these 1's give extensions as desired and the commutative diagram as
desired.

\item This is obvious.

\item Notice that outside $\{A=0\}$, since multiplication by $A$ and
natural inclusions of ideal sheaves are isomorphisms, $\eta_{i+1}$ is
also an isomorphism. So we need to verify the claim at points on
$A=0$. For such a point, which is not on $C$,
$\cI_{C_{i+1}}\hookrightarrow \cI_{C_i}$ is an isomorphism. So the
cokernel of $\eta_{i+1}$ is the same as the cokernel of $\cdot A$, so
claim is proved for such points. Now let $p\in C$. Then near $p$,
$\cI_{C_i}=(z,A^i)$, where $z=x^{pk}$ or $y^{pl}$ at $p$. Also
$\cI_{C_{i+1}}=(z,A^{i+1})$. Pick a basis for $M_i$ and $M_{i+1}$,
which go to $z, A^i, A^{i+1}$. Then $\eta_{i+1}$ is represented by a
matrix of the form $(v_1,v_2)$, where $v_i\in M_i$ and

$$v_1=(1,0)+\lambda\alpha_i(1),\quad v_2=(0,A)+\mu\alpha_i(1),$$

where $\lambda,\mu\in\cO_p$.  But the fact that
$\eta_{i+1}\circ\alpha_{i+1}=\alpha_i\circ A$ implies immediately that
the cokernel of $\eta_{i+1}$ is in fact isomorphic to $\cO_p/A\cO_p$.
(In fact, this argument shows that $M_i/\eta_{i+1}M_{i+1}$ is a line
bundle on the hypersurface $A=0$.  Moreover, one can even write down
exactly this line bundle, though we will not use that fact.) Thus the
map $M_i\otimes \cO_Y(-N)\to M_i$, got by multiplication by $A$,
factors through $\eta_{i+1}$.

\item This follows essentially from the fact that
$H^1(L_i^{-1}\otimes L_{i-1}(-d))=0$. We will construct the $g$'s by
induction. By the stated vanishing, we have $g_1:L_2\to M_1(-d)$ by
lifting $A\in\cI_{C_1}$. Clearly $\phi\circ g_1=\alpha_1\circ A$. So
assume we have constructed $g_{i-1}$ with the required property. So we
have $g_{i-1}\circ A: L_{i+1}\to M_{i-1}(-d)$. This is just the
composite,
$$L_{i+1}=L_i\otimes\cO_Y(-N)\stackrel{g_{i-1}\otimes
1}{\longrightarrow}
M_{i-1}(-d)\otimes\cO_Y(-N)\stackrel{A}{\longrightarrow} M_{i-1}(-d).$$
Now by the previous claim, we see that the last map factors through
$\eta_i$. So we get a map $g_i:L_{i+1}\to M_i(-d)$ such that $\eta_i
g_i=g_{i-1}\circ A$. To compute $\beta_i g_i$ we may compose it with
the natural inclusion of ideals and then it is just
$$\beta_{i-1}\eta_i g_i=\beta_{i-1} g_{i-1}\circ A=A^{i-1}\circ
A=A^i.$$

\item We must show that $f$ is injective at every point. So let $P\in
C$. Then $\alpha_i$'s are all zero. So if $f(x_p,\ldots,x_2)=0$ at
$P$, then $g_i(x_i)=0$. But since near $P$, $A^{i-1}$ is part of a
minimal set of generator of $\cI_{C_{i-1}}$ and thus $g_i(x_i)=0$
implies $x_i=0$ at this point for all $i$. Now let $P\not\in C$. Then
$\alpha_i$'s are injective at this point. If $f(x_p,\ldots, x_2)=0$
at $P$, we will use descending induction to prove that all the
$x_i$'s are zero. Clearly $\alpha_p(x_p)=0$ implies $x_p=0$. Assume
we have proved $x_p=\dots= x_k=0$. Then by looking at the definition
of $f$, we see that $\alpha_{k-1}(x_{k-1})=0$ and thus $x_{k-1}=0$.

\item We should show that ${\rm Im\ }\theta\circ f\subset {\rm Im\ }
f$.
\[
\begin{array}{ll}
\theta\circ f(x_p,\dots,x_2)&\\
=\theta(-\alpha_p(x_p),
-\alpha_{p-1}(x_{p-1})+g_p(x_p), \ldots, -\alpha_2(x_2)+g_3(x_3),
g_2(x_2))&\\
 = f(0,Ax_p,\ldots,Ax_3)&
\end{array}
\]

\item For this it clearly suffices to prove that
$$\cM'={\rm Im\ }f(\cL)+{\rm Im\ }\theta(\cM(-d)) +M_p=\cM.$$
So let $b=(b_p,\dots, b_1)\in\cM$.

First let us look at a point
$P\not\in C$. We will show inductively that there exists a
$c_i\in\cM'$ such that for all $j\geq i$, the $j^{th}$ coordinate of
$b-c_i$ is zero. We may clearly take $c_p=(b_p,0,\ldots,0)$. So by
induction we may assume that $b_j=0$ for $j>i$. Since $P\not\in C$,
we see that $\alpha_i(L_i)+\eta_{i+1}(M_{i+1})=M_i$ at $P$. So we may
write $b_i=\alpha_i(s)+\eta_{i+1}(t)$. Let us first look at the case
when $i\geq 2$. Consider
$$c_i=f(0,\ldots, 0,-s,0,\ldots,0)+\theta(0,\dots,0, t,\dots,0)\in\cM'$$
By our definition of $f,\theta$, we can easily see that $b-c_i$ has
all coordinates upto $i$ zero. Next look at the case when $i=1$.
Again since $P\not\in C$, we see that
$\phi(M_1(-d))+\eta_2(M_2)=M_1$. Thus we can write
$b_1=\phi(s)+\eta_2(t)$. Let $c_1=\theta(0,\dots,0,t,s)\in\cM'$ and
we are done.

Now let us look at points $P\in C$. Now we will show inductively that
there exists $c_i\in\cM'$ such that $b-c_i$ has $j^{th}$ coordinate
zero for all $j\leq i$. For $i=1$, we have
$g_2(L_2(d))+\eta_2(M_2)=M_1$ at $P\in C$. So $b_1=g_2(s)+\eta_2(t)$.
Take $c_1=f(0,\ldots,0,s)+\theta(0,\ldots,0,t,0)\in\cM'$. So assume
that $b_j=0$ for $j<i$. Again let us first look at the case when
$i<p$. Again we have $g_{i+1}(L_{i+1}(1))+\eta_{i+1}(M_{i+1})=M_i$ at
$P\in C$. Therefore we may write $b_i=g_{i+1}(s)+\eta_{i+1}(t)$.
Consider
$c_i=f(0,\ldots,s,\ldots,0)+\theta(0,\ldots,t,\ldots,0)\in\cM'$. One
easily checks that $b-c_i$ has all coordinates zero upto the
$i^{th}$ by using our definition of $f,\theta$. Finally assume $i=p$.
But $(b_p,0,\dots,0)$ clearly belongs to $\cM'$ and thus we are done.

\item If $M_1$ is a direct sum of line bundles, we get $\cI_{C_1}$
is a complete intersection, say of $f,g$ of degrees $a,b$. Then we see
by our Koszul exact sequence, $a+b=d$. Also degree of $C_1=ab$ by
Bezout's theorem. Since $C_1$ is supported along the line $C$, we may
compute its degree by computing the length of $\cO_{C_1}$ at the
generic point of $C$. Easy to see that this is $pkl$.  Thus
$pkl=ab\leq d^2/4$. This contradicts our choice of $N,k,l$.

Thus $\cM$ is not a direct sum of line bundles.  Since $\cL$ is a
direct sum of line bundles, this implies that $M$ is also not a direct
sum of line bundles.

\end{enumerate}

\section{Characteristic zero case}

In this section, we will analyse our construction in characteristic
zero. So assume that our base field is $\C$, the complex numbers from
now on. Let $E$ be a vector bundle on $\P^3$ with a nilpotent
endomorphism $\phi: E\to E(d)$  for a fixed integer $d$ and let $F$ be
the direct sum of two line bundles. Let $\psi:F\to E$ be a
homomorphism such that $\phi(E(-d))+\psi(F)=E$. Assume further that
$E$ is not a direct sum of line bundles and $\rank E=k$.  Assume that
we have chosen $E$ with the smallest possible rank.

\begin{enumerate}
\item {\em We may assume that rank of $\phi$ is $k-1$}

If not rank of $\phi\leq k-2$ and since $F$ has rank two, we see that
$E=\phi(E(-d))\oplus\psi(F).$ Then easy to see that the vector bundle
$\phi(E)$ also has all the properties, with $\phi\psi:F(-d)\to\phi(E)$
replacing $\psi$.  So by minimality of ranks we have proved the claim.
(Note that $E$ is not a direct sum of line bundles implies $\phi(E)$ is
also not a direct sum of line bundles since $F$ is.)

Thus we see that $\ker\phi$ has rank one, $\phi^{k-1}\neq 0$ and
$\phi^k=0$.

\item {\em We may assume $\phi^{k-1}(E)$ is an ideal sheaf defining a
complete intersection curve}

By the previous step, $\phi^{k-1}(E)$ is a rank one torsion free sheaf
and hence isomorphic to $I(l)$ for some ideal sheaf $I$ of height
bigger than or equal to 2 and $l$ an integer. Since we may twist $E$
and $F$ by $\cO(-l)$, we may assume that $l=0$. So we have an exact
sequence,

\begin{equation}\label{basic1}
0\to M\to E\to I\to 0
\end{equation}

where $M=\ker\phi^{k-1}$. Since $\psi(F)$ maps onto $E/\phi(E(-d))$
and since $\phi(E(-d))$ is contained in $M$ we see that $F$ maps onto
$I$. So if $I$ is a proper ideal we would be done. If not $I=\cO$ and
then since $F$ is rank two, this surjection must split. So $F$ is
isomorphic to $\cO\oplus\cO(n)$ for some $n$ and the map $E\to\cO$
also splits.  Thus $E=M\oplus\cO$.  Notice that $\psi\cO(n)\subset M$.
We have

$$M=(\phi(E(-d))+\psi(\cO(n))+\psi(\cO))\cap M=\phi(E(-d))+\psi(\cO(n))$$

which in turn is equal to

$$\phi(M(-d))+\phi(\cO(-d))+\psi(\cO(n)).$$

Thus replacing $(E,\phi,\psi)$ by $M$, $\phi: M(-d)\to M$, the
restriction of $\phi$ and
$$\psi':\cO(-d)\oplus\cO(n)\to M\quad \psi'(x,y)=\phi(x)+\psi(y),$$
we get an example with smaller rank. This contradicts the minimality
of rank assumption.

The fact that $I$ has homological dimension one implies
that $M$ is a vector bundle.

\item {\em  rank $E\geq 3$.}

Clearly rank is bigger than 1, since $E$ is not a direct sum of line
bundles. So if the
assertion is false, then rank must be two. But since $E$ maps onto a
complete intersection ideal, by Serre's construction we see that $E$
must be a direct sum of line bundles.

{\em From now on we assume that the rank of $E$ is three.}

We have the basic diagram,
\begin{equation}
\begin{array}{ccccccccc}\label{basic2}
0&\rightarrow&M&\rarr&E&\rarr&I&\rarr&0\\
&&\theta\uparrow&&\psi\uparrow&&id\uparrow&&\\
0&\rarr&L&\rarr&F&\rarr&I&\rarr&0
\end{array}
\end{equation}

where $L$ is the determinant of $F$. Notice that, since $I$ is a
proper ideal, $\psi$ and hence $\theta$ are both inclusions.

\item {\em $M/\phi(E(-d))$ is not supported along divisors}

Since $\phi$ has rank two, $M/\phi(E(-d))$ is a torsion sheaf. If it
is supported along divisors, choose $Z\subset \supp(M/\phi(E(-d)))$, a
reduced irreducible divisor. Then one sees that
$E/\phi(E(-d))\otimes\cO_Z$ has rank at least two. We of course have a
surjection using $\psi$ from $F\otimes\cO_Z$ to this sheaf. This
implies that this surjection must be an isomorphism since
$E/\phi(E(-d))\otimes\cO_Z$ has rank at least two over $Z$ and
$F_{\mid Z}$ is a vector bundle of rank two. Thus $E\otimes\cO_Z$ is a
direct summ of $F\otimes\cO_Z$ and a line bundle which must be of the
form $\cO_Z(l)$ for some $l$ by determinant considerations. Thus
$E\otimes\cO_Z$ is a direct sum of line bundles. Now one can see
easily that $H^1(E(*))=H^2(E(*))=0$ using the fact $H^1(\cO_Z(*))=0$.
So by Horrocks criterion, we see that $E$ is also a direct sum of line
bundles leading to a contradiction.

Thus we see that $\det E=\det M=\det\phi(E(-d))$. But the kernel of
$\phi:E(-d)\to E$ is a rank one sheaf which is reflexive and hence a
line bundle, say $A$. So $\det E(-d)=A\otimes\det\phi(E(-d))$ and
hence $A=\cO(-3d)$.
Notice that $A\subset M(-d)$ and $M(-d)/A\cong J(l)$ for some ideal
sheaf $J$
of height bigger than one and $l$ an integer. Also $\phi$ restricts
to a nilpotent endomorphism $M\to M(d)$ and one easily sees that it is
obtained by going to $J(l)$ and composing it with some embedding of
$J(l)$ in $A(d)$. In particular, $l\leq -2d$. On the other hand we have a
natural map $I=E/M\to M(d)/A(2d)=J(l+2d)$ induced by $\phi$ and this is
injective. So $l+2d\geq 0$ and thus $l=-2d$. Thus $\det M=\cO(-3d)$ and
hence $\det E=\det M=\cO(-3d)$.

So we have the next basic exact sequence,
\begin{equation}\label{basic3}
0\rarr\cO(-d)\rarr M(d)\rarr J\rarr 0
\end{equation}

\item {\em $J$ is a proper local complete intersection ideal of a curve}

We only need to show that $J$ is a proper ideal since $M$ is a rank
two vector bundle. If $J=\cO$, then $M(d)=\cO(-d)\oplus\cO$.
The image of $\phi(E)$ in $J$ is $I$. We have
$$M(d)= E(d)\cap M(d)=(\phi(E)+\psi(F(d)))\cap M(d)$$
$$=\phi(E)+\psi(F(d))\cap M(d)=\phi(E)+\theta(L(d))$$

In particular, $J=I+{\rm image}(L(d))$. If $J=\cO$ then, we see that
$F\oplus L(d)$ surjects onto $J$ and thus one of them must be $\cO$
(Any three non-trivial polynomials in $\P^3$ have a common zero). This
is impossible since $I$ is non-trivial and $L$ is the determinant of
$F$ unless $d=a+b$ where $F=\cO(-a)\oplus\cO(-b)$. This implies in
particular that $d>0$ and going to the commutative diagram
(\ref{basic2}) above, using the fact that $M=\cO(-d)\oplus\cO(-2d)$
and $\theta$ is injective, we get an exact sequence,

$$0\to F\stackrel{\psi}{\longrightarrow} E\to \cO(-2d)\to 0,$$

which implies that $E$ is a direct sum of line bundles.  We see  $J$
is a proper ideal of height at least two and since $M$ surjects onto
it, it must be a local complete intersection ideal of height two. The
image of $L(d)$ corresponds to an element $G\in J$ of degree
$c=a+b-d.$ We have $J=I+(G)$. We have an exact sequence,

\begin{equation}\label{basic5}
0\rarr I\rarr J\rarr T\rarr 0\label{basic4}
\end{equation}

where $T=\cO_X(-c)$, with $X=\supp T$. Let us also denote
by $C$, the curve defined by $I$ and $D$ that defined by $J$. Then $D$
is a closed subscheme of $C$.

\item {\em $X$ is a local complete intersection curve}

First we show that it has no components of dimension two. If it
did, say $Z$, an irreducible hypersurface, restricting to $Z$, we see
that $M\otimes\cO_Z$ surjects onto $\cO_Z(-c)$ and then
$M\otimes\cO_Z$ must be isomorphic to $\cO_Z(-c)\oplus\cO_Z(c)$.
Again as before we get by Horrock's criterion, that
$M\cong\cO(c)\oplus\cO(-c)$ and since it is supposed to surject onto
$J$ and $c\geq 0$ or $-c\geq 0$, we reach a contradiction. To check
that it is a local complete intersection curve, let us do it locally.
If $x\notin C$ then since $D\subset C$, we see that $J$ and $I$ are
locally $\cO$ and if the inclusion were not an isomorphism, the
cokernel $T$ would be supported on a hypersurface which we have seen
is not the case. So $X_{\rm red}\subset C$. If $x\in C$ and not in $D$
clearly we are done since $I$ is a local complete intersection. So let
us look at a point $x\in D$. Since $I$ and $J$ are local complete
intersections and $J/I$ is principal, we may choose generators so that
$J=(f,g)$ and $I=(f,gh)$ and then $X$ is defined by $(f,h)$ and we are
done.

Dualising the above exact sequence (\ref{basic5}), we get,
$$ \sext^1(T, \cO)\rarr\sext^1(J,\cO)\rarr\sext^1(I,\cO)\rarr
\sext^2(T, \cO)\rarr 0$$
Since $X$ is a local complete intersection, we get that the first
term is zero. Also $\sext^2(T, \cO)\cong\omega_X(c+4)$.
{}From our basic exact sequences \ref{basic1}, \ref{basic3} it is easy
to see that
$\sext^1(I,\cO)\cong \cO_C(a+b)$ and $\sext^1(J,\cO)\cong \cO_D(d)$.
Consider the diagram,
\[
\begin{array}{ccccccccc}
0&\rarr&\cO(-(a+b))&\rarr&F&\rarr&I&\rarr&0\\
&&G\downarrow&&\phi\psi\downarrow&&\downarrow&&\\
0&\rarr&\cO(-d)&\rarr&M(d)&\rarr&J&\rarr&0
\end{array}
\]

This diagram is commutative. The only fact we need to show is that the
vertical arrow on the left is multiplication by $G$. But the map $\psi$
restricted to $L=\cO(-(a+b))$ maps to $M$ and its image in $J$ is $G$.
Since $\phi:M\to M(d)$ is got by going to $J(-d)$ and including it in
the kernel (of $\phi: M(d)\to M(2d)$), $\cO(-d)$, we see that
$\phi\psi$ restricted to $L$ is precisely multiplication by $G$.
The fact that the vertical map on the left is multiplication by $G$
implies by dualising, that the natural map $\cO_D(d)\to \cO_C(a+b)$ is
multiplication by $G$. So the cokernel is isomorphic to $\cO_D(a+b)$,
since the image of $G$ generates $J$ in $\cO_C$. So we get that
$\omega_X(c+4)\cong\cO_D(a+b)$. In particular $X=D$. So $C$ is
set-theoretically the same as $D$. Thus $G^2\in I$. Let $I=(f,g)$.
Write $G^2=Af+Bg$.

\item {\em $A, B, f, g$ have no common zeroes in $\P^3$}

If $p$ is such a zero, since $f,g$ vanish there, $p\in D$. Locally at
that point, we see that $J$ is generated by one of $f,g$ and $G$. Let
us assume without loss of generality, that $(f,G)=J$. Then $I=(f,
G^2)$. So $(f,g)=(f, Af+Bg)=(f,Bg)$ and thus $B$ must be a unit
there. That means $B$ does not vanish at $p$.

\item {\em The ideal sheaf defined by $f,g$ in $G^2=0$ is a line bundle}

This is obvious by the earlier remark and a local checking.

\begin{lemma}\cite{Pesk}
Let $Y\subset\P^3$ be the hypersurface defined by $G=0$ and $Y'$ that
defined by $G^2=0$. Then any line bundle on $Y'$ is trivial. In other
words, the natural map $\Pic \P^3\to\Pic Y'$ is an isomorphism.
\end{lemma}

\begin{proof} Proof is essentially due to Ellingsrud {\em et. al.} One
can first reduce to the case when $Y$ is a reduced irreducible
hypersurface as follows.  First we may assume that $Y$ is irreducible
(not necessarily reduced).  For this let $Y_1, Y_2,\ldots ,Y_n$ be the
scheme theoretic irreducible components of $Y'$. Let $L\in\Pic Y'$
such that its restriction to all the $Y_i$'s are of the form
$\cO_{Y_i}(n_i)$ for some $n_i$. Then restricting to the intersection
$Y_i\cap Y_j$ we see that $n_i$ is a constant independent of $i$.
Twisting by the negative of this line bundle, we are reduced to
proving that if $L$ restricted to all $Y_i$ is $\cO_{Y_i}$, then
$L=\cO_{Y'}$. The proof is standard boot strapping.  Assume we have
proved that $L$ restricted to $Z=Y_1\cup Y_2\cup\ldots Y_k$ is trivial
($\cong\cO_Z$). We want to show that $L$ restricted to $Z'=Z\cup
Y_{k+1}$ is also trivial. We have the natural exact sequence,

$$
0\to \cO_{Y_{k+1}}(-Z)\to\cO_{Z'}\to\cO_Z\to 0
$$

Tensoring this by $L$ and noting that $L_{\mid Z}\cong\cO_Z$ and $L$
restricted to $Y_{k+1}$ is trivial and hence the kernel has no first
cohomology, we see that the section $`1'\in H^0(\cO_Z)$ can be lifted
to a section of $L$ restricted to $Z'$. Easy to see that this section
generates this line bundle by using the fact that $Z\cap
Y_{k+1}\neq\emptyset$.

So let us assume that $Y$ is irreducible defined by $G=0$. By the
$\cO^{*}$ exact sequences one can see that $\Pic Y'\to\Pic Y$ is
injective. Thus we reduce to the case of $Y$ a reduced irreducible
surface. Notice that $\Pic Y/(\Z=\Pic \P^3)$ is torsion free. Since
$H^1(\cO_{Y})=0$, it follows that the natural map

$$\Pic Y=H^1(\cO_Y^*)\stackrel{dlog}{\longrightarrow}
H^1(\Omega^1_{Y})$$

is injective [see \cite{SGA6}, Th\'{e}or\`{e}me 4.7, ({\romannumeral
3})].  So we only need to understand the map

$$
H^1(\Omega^1_{Y'}\otimes\cO_Y)\to H^1(\Omega^1_Y).
$$

We have the fundamental exact sequence,
$$ \cO_{\P^3}(-Y')\otimes\cO_{Y'}\to \Omega^1_{\P^3}\otimes\cO_{Y'}\to
\Omega^1_{Y'}\to 0$$

When we restrict this to $Y$, we see that the first map is zero, since
it is given essentially by the derivative of $G^2$, which is zero
modulo $G$.

Thus,
$\Omega^1_{Y'}\otimes\cO_Y\cong\Omega^1_{\P^3}\otimes\cO_Y.$ So we
have,

\begin{equation}\label{one}
H^1(\Omega^1_{Y'}\otimes\cO_Y)\cong
H^1(\Omega^1_{\P^3}\otimes\cO_Y)\cong \C
\end{equation}

We have the standard commutative diagram,

\[
\begin{array}{ccccc}
\Pic\P^3&\into&\Pic Y'&\into&\Pic Y\\
\downarrow&&\downarrow&&\downarrow\\
H^1(\Omega^1_{\P^3})&\to&H^1(\Omega^1_{Y'})&\to&H^1(\Omega^1_Y)
\end{array}
\]
{}From the various inclusions, we see that if $\alpha\in\Pic Y'$, we
must show that its image in $H^1(\Omega_Y)$ comes from $\Pic\P^3$. By
chasing the diagram, and using the  equation (\ref{one}) we see that
there exists a {\em complex number},
$z$, such that $\alpha=zH$ where $H$ is the tautological class in
$\P^3$. By restricting to a general map from a smooth curve to $Y'$
and noting that then the classes correspond to degrees of the
corresponding line bundles on this curve, we get
$z\in\Q$. But now using the torsion freeness of $\Pic Y/\Z$, we see
that $\alpha\in\Pic \P^3$.  Thus we get the result.
\end{proof}

Putting all these together, we see that $I=(f,g)=(G^2,h)$ for some
$h$. But then $J=(f, g,G)=(G,h)$ is a complete intersection.  Let
$\deg h=e$ (where $e=a$ or $b$).  $M(d)=\cO(-e)\oplus \cO(-c)$. Then
let us look at the map $\theta$.  It is a map from $\cO(-a-b)\to
\cO(-e-d)+\cO(-c-d)$. But $-c-d=-a-b$ and thus the quotient must be
either a line bundle or must be a line bundle direct sum a sheaf $\cF$
supported along a divisor. In the former case, since we also know
that it is the quotient of $E$ by $\psi(F)$, we get $E$ to be a direct
sum of line bundles, which we have assumed is not the case. In the
latter case one can see that $E$ can not be a vector bundle at any
point of intersection of $\supp \cF$ and the curve defined by $f,g$.
This is also a contradiction. This finishes the proof that such
examples can not exist.
\end{enumerate}

Thus in characteristic zero, any such $E$ must have rank at least 4.
There are several technical difficulties which must be overcome to
analyse these cases (say of rank 4). The analysis which we have
carried out will be discussed elsewhere.

\bibliography{pn}
\bibliographystyle{alpha}
\end{document}